\documentclass[12pt]{amsart}
\usepackage[centertags]{amsmath}

\usepackage{amsfonts}
\usepackage{amssymb}
\def\pa{{\sf PA}}
\def\MM{{\mathcal M}}

\def\res{\hspace{-4pt} \upharpoonright \hspace{-2pt}}

   \DeclareMathOperator{\Th}{Th}
      \DeclareMathOperator{\aut}{Aut}
      \def\fA{\mathfrak A}     \def\fB{\mathfrak B}

      \def\into{\longrightarrow}

\title{The Automorphism Group of \\ a Resplendent Model}
\author{James H. Schmerl }
\date{\today}

\begin{document}

\maketitle

Our purpose here is to prove the following theorem.

\bigskip

{\sc Theorem} 1: {\em If ${\fA}$ is an infinite resplendent structure for some 
some finite first-order language, then $\Th(\aut(\fA))$ is undecidable.}

\bigskip

To clarify notation, $\aut(\fA)$ is the group of all automorphisms of $\fA$
and $\Th(\aut(\fA))$ is the first-order theory of that group. An ${\mathcal L}$-structure $\fA = (A, \ldots)$, where ${\mathcal L}$ is  some first-order language, is {\it resplendent}
(consult, for example,  \cite{hod} or \cite{poi}) if, whenever ${\mathcal L}' \supseteq {\mathcal L}$ is a language that may have some constant symbols  denoting 
 elements of $A$ and $\sigma$ is an 
${\mathcal L}'$-sentence that has a model ${\fB}'$ such that ${\fA} \preccurlyeq 
{\fB}' \res {\mathcal L}$, then ${\fA}$ has an ${\mathcal L}'$-expansion ${\fA}'$ that is a model of~$\sigma$. 

While resplendent structures are perhaps not that well known (although \cite{kos} is an attempt to change that), the countable 
ones are more so, being precisely the countable, recursively saturated structures. 
It should be pointed out that every countable first-order theory $T$ having an infinite model has  resplendent models of every infinite cardinality $\kappa \geq |T|$.

The impetus for proving this theorem comes from two sources. The first is  the remarkable theorem of 
Bludov, Giraudet, Glass and Sabbagh \cite{bggs} that states: {\em If $T$ is an 
arbitrary first-order theory having an infinite model, then it has a model 
$\fA$ such that $\Th(\aut(\fA))$ is undecidable.} In this result, the language 
is not restricted to being finite. Indeed, this result has the following 
addendum: {\em Moreover, if $\kappa \geq |T|$ is an infinite cardinal, then 
$\fA$ can be chosen to have cardinality $\kappa$.} This theorem suggests 
the problem of finding a nice class of models  ${\fA} \models T$ for which 
$\Th(\aut(\fA))$ is undecidable. 

The second source is the following 
theorem that I proved (Theorem~6.1 of \cite{s}) about models of Peano Arithmetic (\pa): {\em If $\MM$ is a countable, recursively saturated model of 
\pa, then $\Th(\aut(\MM))$ is undecidable.}  Although the proof of this result 
was very specific to models of \pa, it, together with the  theorem from \cite{bggs},
suggests that \pa\ might have nothing to do with the result. Indeed it doesn't,
as  Theorem~1 makes clear.

The proof of  Theorem~1 relies heavily on the proof in \cite{bggs}. 
The automorphism group of the model 
 of $T$  obtained in the proof in \cite{bggs} had an undecidable universal Horn theory.  Similarly,  the proof of Theorem~1 will show that the universal Horn theory  of $\aut(\fA)$ is undecidable.
The crucial lemma from \cite{bggs} states that there is a 
finitely presented group $G$ with an unsolvable word problem that is 
embeddable in $\aut(({\mathbb Q},<))$. (Here, $({\mathbb Q},<)$ is the rationals with their usual ordering.) As pointed out in \cite{bggs}, if 
$H$ is any group that has $G$ as a subgroup, then $\Th(H)$ is undecidable.
Thus, to prove  Theorem~1 just for  countable ${\fA}$, it suffices to prove the following 
theorem.

\bigskip

{\sc Theorem 2}: {\em If ${\fA}$ is a countably infinite, recursively saturated 
structure for a finite language, then $\aut(({\mathbb Q},<))$ is embeddable 
in $\aut(\fA)$.}

\bigskip

{\it Proof}. We give a proof using models of \pa. For a countable, recursively saturated model $\MM = (M,+,\times,0,1,\leq)$ of \pa, it is well known that 
$\aut(({\mathbb Q},<))$ is embeddable 
in $\aut(\MM)$ and even in $\aut((\MM,a))$ for any $a \in M$.  (See, for example, Corollary~5.5.2 of \cite{ksbook}.)
Now, by the resplendency of $\fA$ and the Arithmetized Completeness Theorem (which is provable in \pa) there is a recursively saturated model 
$\MM \models \pa$ in which $\fA$ is  $\Delta_2$ and $A = M$. 
Thus, there is some $a \in M$ from which $\fA$ is definable in $\MM$. 
Then, since $\aut((\MM,a))$ is a subgroup of $\aut(\fA)$, we get that 
$\aut(({\mathbb Q},<))$ is embeddable in $\aut(\fA)$.  \qed

\bigskip

{\it Proof of Theorem~1}. Let 
$$G = \langle g_1,g_2, \ldots, g_m \ |\ 
u_1(\overline g) = 1, \ldots, u_n(\overline g) = 1 \rangle$$ be a presentation 
of $G$. If  $w = w(\overline g)$ is a word, then let $\sigma_w$ be the 
universal Horn sentence 
$$
\forall x_1,x_2, \ldots, x_m[\big( \bigwedge_{i=1}^n u_i(x_1,x_2, \ldots, x_m) = 1\big) \into w(x_1,x_2, \ldots, x_m) = 1]
$$
in the language of group theory. The undecidability of $\Th(\aut(\fA))$ will now follow from the following claim.

\smallskip

{\it Claim}: For any word $w = w(\overline g)$,
$\aut(\fA) \models \sigma_w$ iff $G \models w(\overline g) = 1$. 

\smallskip

First,
suppose that  $G \models w(\overline g) = 1$. Let 
$f_1,f_2, \ldots,f_m \in \aut(\fA)$ be such that   $\aut(\fA)\models  u_i( f_i) = 1$ 
for each $i$. Consider the homomorphism from $G$ into $\aut(\fA)$ that maps $g_i$ to $f_i$. This homomorphism maps $w(\overline g)$ to $w(\overline f)$, so that $w(\overline f) =1$. Therefore, $\aut(\fA) \models \sigma_w$

Conversely, suppose that $G \models w(\overline g) \neq 1$.  
Let ${\mathfrak B} \equiv {\fA}$ be countable and recursively saturated, 
so by Theorem~2, $G$ is embeddable in $\aut(\mathfrak B)$. 
Thus, by the resplendency of $\fA$,
there are 
$f_1,f_2, \ldots,f_m \in \aut(\fA)$ such that $\aut(\fA)\models w(\overline f) \neq 1 \wedge \bigwedge_{i=1}^n u_i(\overline f) = 1$. Clearly, then $\aut(\fA) \models \neg\sigma_w$, proving the claim. \qed
 

\bibliographystyle{plain}

\end{document}